\documentclass[10pt]{article}

\usepackage{amsmath, amssymb, amscd}

\def\eqref#1{(\ref{#1})}

\newcommand{\arrow}{{\:\longrightarrow\:}}
\newcommand{\Z}{{\Bbb Z}}
\newcommand{\C}{{\Bbb C}}
\newcommand{\R}{{\Bbb R}}

\newcommand{\6}{\partial}

\newcommand{\restrict}[1]{{\left|_{{\phantom{|}\!\!}_{#1}}\right.}}

\renewcommand{\c}[1]{{\cal #1}}
\newcommand{\calo}{{\cal O}}

\renewcommand{\tilde}{\widetilde}
\renewcommand{\bar}{\overline}
\renewcommand{\phi}{\varphi}
\renewcommand{\epsilon}{\varepsilon}
\renewcommand{\geq}{\geqslant}
\renewcommand{\leq}{\leqslant}

\newcommand{\fl}{{\rm fl}}

\newcommand{\End}{\operatorname{End}}
\newcommand{\Tot}{\operatorname{Tot}}

\newcommand{\Hom}{\operatorname{Hom}}
\newcommand{\Coh}{\operatorname{Coh}}
\newcommand{\Sup}{\operatorname{Sup}}
\newcommand{\Ext}{\operatorname{Ext}}

\newcommand{\Sec}{\operatorname{Sec}}

\newcommand{\Sing}{\operatorname{Sing}}
\newcommand{\codim}{\operatorname{codim}}

\newcommand{\rk}{\operatorname{rk}}

\newcommand{\Def}{\operatorname{Def}}
\newcommand{\Tw}{\operatorname{Tw}}
\newcommand{\Tr}{\operatorname{Tr}}
\newcommand{\HolBun}{\operatorname{HolBun}}
\newcommand{\Bun}{\operatorname{Bun}}
\newcommand{\Fl}{\operatorname{Fl}}
\newcommand{\Spec}{\operatorname{Spec}}

\renewcommand{\Re}{\operatorname{Re}}

\newcommand{\comment}[1]{{}}

\def\blacksquare{\hbox{\vrule width 4pt height 4pt depth 0pt}}
\def\endproof{\blacksquare}

\makeatletter

\@ifundefined{Bbb}
     {\newcommand{\Bbb}[1]{{\mathbb #1}}}%
{}

\newcommand{\ps@verbit}{%
  \renewcommand{\@oddhead}{%
          \scriptsize
          {Coherent sheaves on generic tori}
          \hfil\tiny {M. Verbitsky, \ \ \ \ \ October 19, 2003 }}
  \renewcommand{\@evenhead}{\@oddhead}
  \renewcommand{\@oddfoot}{\hfil\thepage\hfil}
  \renewcommand{\@evenfoot}{\@oddfoot}}
 
\newcounter{Mycounter}[section]
\newcounter{lemma}[section]
\renewcommand{\thelemma}{{Lemma \thesection.\arabic{lemma}}}
\newcommand{\lemma}{%
     \setcounter{lemma}{\value{Mycounter}}
     \refstepcounter{lemma}
     \stepcounter{Mycounter}
     {\bf \thelemma:\ }}

\newcounter{claim}[section]
\renewcommand{\theclaim}{{Claim \thesection.\arabic{claim}}}
\newcommand{\claim}{%
     \setcounter{claim}{\value{Mycounter}}
     \refstepcounter{claim}
     \stepcounter{Mycounter}
     {\bf \theclaim:\ }}

\newcounter{sublemma}[section]
\newcounter{corollary}[section]
\renewcommand{\thecorollary}{{Corollary \thesection.\arabic{corollary}}}
\newcommand{\corollary}{%
     \setcounter{corollary}{\value{Mycounter}}
     \refstepcounter{corollary}
     \stepcounter{Mycounter}
     {\bf \thecorollary:\ }}

\newcounter{theorem}[section]
\renewcommand{\thetheorem}{{Theorem \thesection.\arabic{theorem}}}
\newcommand{\theorem}{%
     \setcounter{theorem}{\value{Mycounter}}
     \refstepcounter{theorem}
     \stepcounter{Mycounter}
     {\bf \thetheorem:\ }}

\newcounter{conjecture}[section]
\newcounter{proposition}[section]
\renewcommand{\theproposition}
       {{Proposition \thesection.\arabic{proposition}}}
\newcommand{\proposition}{%
     \setcounter{proposition}{\value{Mycounter}}
     \refstepcounter{proposition}
     \stepcounter{Mycounter}
     {\bf \theproposition:\ }}

\newcounter{definition}[section]
\renewcommand{\thedefinition}
       {{Definition~\thesection.\arabic{definition}}}
\newcommand{\definition}{%
     \setcounter{definition}{\value{Mycounter}}
     \refstepcounter{definition}
     \stepcounter{Mycounter}
     {\bf \thedefinition:\ }}

\newcounter{example}[section]
\newcounter{remark}[section]
\renewcommand{\theremark}{{Remark \thesection.\arabic{remark}}}
\newcommand{\remark}{%
     \setcounter{remark}{\value{Mycounter}}
     \refstepcounter{remark}
     \stepcounter{Mycounter}
     {\bf \theremark:\ }}

\newcounter{problem}[section]
\newcounter{question}[section]
\@addtoreset{equation}{section}
\@addtoreset{footnote}{section}
\makeatother

\begin{document}

\begin{center}
{\LARGE\bf
Coherent sheaves on generic compact tori
}
\\[4mm]
Misha Verbitsky,\footnote{ The author is EPSRC advanced fellow 
supported by CRDF grant RM1-2354-MO02 and EPSRC grant  GR/R77773/01
}\\[4mm]
{\tt  verbit@maths.gla.ac.uk, \ \  verbit@mccme.ru}
\end{center}


{\small 
\hspace{0.15\linewidth}
\begin{minipage}[t]{0.7\linewidth}
{\bf Abstract} \\
Let $T$ be a compact complex torus, $\dim T>2$. 
We show that the category of coherent sheaves
on $T$ is independent of the choice of the
complex structure, if this complex structure is generic.
The proof is independent of math.AG/0205210,
where the same result was proven for K3 surfaces
and even-dimensional tori. 
\end{minipage}
}

{
\small
\tableofcontents
}


\section{Introduction}
\label{_Intro_Section_}


Let $M$ be a compact K\"ahler manifold. Category 
$\Coh(M)$ of coherent sheaves on $M$ is an important
invariant $M$ which has deep physical meaning.
For manifolds with ample or anti-ample canonical
class, $M$ can be reconstructed from the
corresponding derived category $D\Coh(M)$
(\cite{_Bondal_Orlov_}, \cite{_Posic_}).
In physics, the objects of $D\Coh(M)$ are interpreted
as certain branes on $M$. Under Mirror
Symmetry, this category corresponds to 
a triangulated category which
is constructed by Fukaya from the symplectic geometry
of Mirror dual manifold $\hat M$. 

Physicists believe that a Mirror dual of $n$-dimensional
compact complex torus is again an $n$-dimensional
compact complex torus. This is one of the first
cases where Mirror Symmetry was known explicitly. 

In \cite{_Verbitsky:cohe_}, we have studied the category
$\Coh(M)$ for $M$ a K3 surface or even-dimensional compact
torus. We proved that $\Coh(M)$ is independent from $M$,
provided that $M$ is generic in its deformation class.
The proof of this statement involves hyperk\"ahler geometry.

Given a hyperk\"ahler structure on $M$, we consider the
corresponding twistor space $\Tw(M)$. To every object 
$F\in \Coh(M)$, we associate a coherent sheaf
$\hat F$ on $\Tw(M)$, in such a way that
the restriction of $\hat F$ to $M$, identified 
with a fiber of the twistor projection
$\pi:\; \Tw(M) \arrow \C P^1$, gives $F$.
Restricting $\tilde F$ to another fiber close
to $M$, we obtain a coherent sheaf on a
manifold which is deformationally equivalent
to $M$. We are free in the choice of
hyperk\"ahler structure, and this allows us
to obtain any given deformation $M'$ of $M$
after several iterations of this procedure. 
This construction gives a functor from
$\Coh(M)$ to $\Coh(M')$. This functor
is equivalence if $M'$ is also generic.

This proof is not very satisfactory, for
several reasons. The equivalence 
$\Coh(M)\cong \Coh(M')$ is not canonical,
because it involves many intermediate choices
(such as the choice of a hyperk\"ahler structure,
and intermediate hyperk\"ahler structures
leading from $M$ to $M'$ as we indicated). 
This dependence is very difficult to estimate.
Also, several important steps of the proof
do not involve hyperk\"ahler geometry at all,
giving an impression that hyperk\"ahler geometry
is mostly irrelevant to the big picture.
Finally, the isomorphism $\Coh(M)\cong \Coh(M')$ 
in \cite{_Verbitsky:cohe_}
is proven only for even-dimensional tori, because
odd-dimensional tori are not hyperk\"ahler.
By avoiding the hyperk\"ahler geometry entirely,
we are able to produce an isomorphism
$\Coh(M)\cong \Coh(M')$ for generic tori of 
arbitrary dimension $>2$. 

\hfill

We use the following notion of generic.

\hfill

\definition\label{_gene_tori_Definition_}
Let $T$ be a compact complex torus, $\dim T \geq 3$.
We say that $T$ is generic, if $T$ has no non-trivial
subtori, 
$H^{1,1}(T) \cap H^{2}(T, \Z) =0$, and
$H^{2,2}(T) \cap H^{4}(T, \Z)=0$. 
As \ref{_non-generic_Proposition_} 
indicates, the set of non-generic tori
is a set of measure zero.

\hfill

In this paper, we prove the following theorem.

\hfill

\theorem\label{_main_theorem_intro_Theorem_}
Let $T$, $T'$ be compact complex tori of dimension
$n\geq 3$. Assume that $T$, $T'$ are generic.
Then the corresponding categories of coherent sheaves
are equivalent: $\Coh(T)\cong \Coh(T')$.

\hfill

We prove \ref{_main_theorem_intro_Theorem_} in 
Section \ref{_moduli_conne_Section_}.

\hfill

Notice that in this paper all  tori 
are assumed to have $\dim T \geq 3$.
We have dealt with the case of 2-dimensional tori in
\cite{_Verbitsky:cohe_}.

\hfill

The proof of \ref{_main_theorem_intro_Theorem_} is based on a concept
of ``$SO(2n)/U(n)$-twistor space'' for a torus. Given 
a compact torus $T$, with flat K\"ahler metric $g$,
consider the space $S$ of all complex structures
$I$ such that $g$ is K\"ahler on $(T,I)$.
Clearly, $S$ is isomorphic to a K\"ahler
symmetric space $SO(2n)/SU(n)$ (so-called
isotropic Grassmanian space). The natural
twistor fibration \[ \Tw(T) \arrow S\] is holomorphic
(\ref{_twi_inte_Proposition_}), and in many ways analogous to the twistor
fibration known from the hyperk\"ahler geometry. This
analogy can be used to adapt the proof 
of equivalence $\Coh(M)\cong \Coh(M')$
{}from \cite{_Verbitsky:cohe_} to the 
present needs. 

This paper is independent from \cite{_Verbitsky:cohe_},
and can be used as introductory reading on the subject,
as the proofs are significantly simplified.

\hfill

We also prove the following results about generic tori.

\hfill

\theorem\label{_geometry_torus_intro_Theorem_}
Let $T$ be a generic torus\footnote{Please notice that 
everywhere in this paper we assume that $\dim T \geq 3$.}.
Then 
\begin{description}
\item[(i)] All proper complex subvarieties of $T$ have dimension 0.
\item[(ii)] Any holomorphic 
vector bunlde on $T$ admits a flat connection compatible with
the holomorphic structure.
\item[(iii)] For any coherent sheaf $F$ on $T$, the reflexization
\[
 F^{**}:= \Hom(\Hom(F, \calo_T), \calo_T)
\]
is a vector bundle.
\end{description}

{\bf Proof:}  \ref{_geometry_torus_intro_Theorem_} (i) is
\ref{_subva_trivi_Proposition_}, \ref{_geometry_torus_intro_Theorem_} (ii) 
is \ref{_bun_admits_flat_Corollary_}, 
and \ref{_geometry_torus_intro_Theorem_} (iii)
is \ref{_gene_refle_Proposition_}. \endproof

\hfill

We shall use the following bit of 
notation, which is due to P. Deligne.
Consider a complex vector space $V$. 
We consider $V$ as a real vector space $V_\R$ with
an action of $\C^*=U(1)\times \R^{>0}$.
Consider some tensor power $W$ of $V_\R$
(for instance, the space $\Lambda^i(V_\R)$
of real exterior differential forms), and
let $W_\C$ be its complexification.
The group $\C^*$ acts 
on $W$ by multiplicativity:
\begin{equation}\label{_Hodge_action_Equation_}
\rho_W:\; \C^*\times W \arrow  W.
\end{equation}
 We say that
a vector $v\in W_\C$ has Hodge type
$(p,q)$ if $\rho(r, v) = r^{p+q}v$,
and $ \rho(\theta, v) = \theta^{p-q}v$,
where $r\in \R^{>0}$ and $\theta\in U(1)$
are standard generators in $U(1)\times \R^{>0}=\C^*$.

This definition is compatible with the
standard notion of Hodge decomposition.


\section{Generic tori}
\label{_gene_tori_Section_}


In this Section, we prove that \ref{_gene_tori_Definition_}
is consistent, that is, the set of non-generic tori
has measure zero in the corresponding moduli space.

Let $R$ be the space of group homomorphisms 
$\Z^{2n}\stackrel \phi\hookrightarrow \C^n$ such that 
$\C^n/\phi(\Z^{2n})$ is compact, up to 
$GL(n, \C)$-action. The space $R$ admits
a natural open embedding to $(\C^n)^{2n}/GL(n, \C)$,
\begin{equation}\label{_embe_to_C^n^2n_Equation_}
R\arrow (\C^n)^{2n}, \ \ \  \phi \arrow (\phi(\alpha_1),
\phi(\alpha_2), ...,  \phi(\alpha_{2n})),
\end{equation}
where $\alpha_i$ are generators of $\Z^{2n}$. 
Clearly, $R$ is identified with the moduli space of 
pairs $(T, \alpha_1, \alpha_2, ... \alpha_{2n})$,
where $T$ is a torus, and $\alpha_1, \alpha_2, ..., \alpha_{2n}$
a basis in $H^1(T, \Z)$. The space $R$ is called
{\bf the marked moduli space of complex tori}.
It is a covering over the usual moduli of
complex tori, with monodromy $GL(2n, \Z)$.

\hfill

\proposition \label{_non-generic_Proposition_}
Let $R_{sp}\subset R$ be the set of all tori which are not
generic, in the sense of \ref{_gene_tori_Definition_}.
Then $R_{sp}$ is a countable union of closed complex
subvarieties of positive codimension in $R$.

\hfill

{\bf Proof:} Given a proper sublattice $L\subset \Z^{2n}$, 
$\rk L = 2l$, let $Z_L\subset R$ be the set of all homomorphisms 
$\Z^{2n}\stackrel \phi\hookrightarrow \C^n$
such that $\phi(L)$ lies inside a complex
space of rank $l$. For any $\alpha\in \Lambda^2_\Z(\Z^{2n})$,
or in $\Lambda^4_\Z(\Z^{2n})$, denote by $Z_\alpha$ the set of all
homomorphisms $\Z^{2n}\stackrel \phi\hookrightarrow \C^n$
such that the element corresponding to $\alpha$
in the cohomology of the corresponding torus $T$ 
is of Hodge type $(1,1)$ or $(2,2)$. Clearly,
$R_{sp}= \bigcup Z_L \bigcup Z_\alpha$, where 
the union is taken over all $L$ and $\alpha$.
Clearly, $Z_L$, $Z_\alpha$ are closed complex
subvarieties of $R$. To prove 
\ref{_non-generic_Proposition_}, it remains
to show that $\codim Z_L>0, \codim Z_\alpha>0$,
for all $L$ and $\alpha$. Since $R$ 
is irreducible, we have $\codim Z_L>0, \codim Z_\alpha>0$
unless $Z_L, Z_\alpha=R$ for some $L$ or $\alpha$.
Given a proper lattice $L\subset \Z^{2n}$, $\rk L = 2l$,
it is easy to find a fomomorphism  
$\Z^{2n}\stackrel \phi\hookrightarrow \C^n$
such that $\C \cdot \phi(\alpha)$ generates a subspace of
dimension $\min(2l, n)$. Therefore, $Z_L\neq R$,
and $Z_L$ has positive codimension. To 
prove that $Z_\alpha$ has positive codimension,
we need to find a torus $T$ with $\alpha\notin H^{p,p}(T)$
for a given $\alpha$. We use the following (more general) lemma.

\hfill

\lemma \label{_gene_torus_no_Hodge_Lemma_}
For any $n$, there exists an $n$-dimensional
compact complex torus $T$ which satisfies
\begin{equation}\label{_very_gene_Equation_}
H^{p,p}(T)\cap H^{2p}(T, \Z)=0 \text{\ \ for all\ \ } 0<p<n.
\end{equation}

\hfill

{\bf Proof:} We need to show that 
for each $v\in H^{2p}(T, \R)$, there exists
a complex structure $I$ on $T$ for which
$v \notin H^{2p}(T_I, \Z)$. Suppose that
for some $v$ this is not true.
Let $\phi\in R$ be a homomorphism
$\Z^{2n}\stackrel \phi\hookrightarrow \C^n=V$, and
$\rho$ the $U(1)$-action on 
$\Lambda^{2p}_\R V$ corresponding to the Hodge
decomposition as in \eqref{_Hodge_action_Equation_}.
Consider the group $G$ generated by such $U(1)$-actions.
It is easy to see that $G=SL(V,\R)$. Then,
for $v$ as above, $v$ is $SL(V, \R)$-invariant.
As follows from classical invariants theory,
any non-trivial $SL(V, \R)$-invariant form on $V$ is 
proportional to the volume form.
This proves \ref {_gene_torus_no_Hodge_Lemma_}.
\endproof

\hfill

\remark
As the proof of \ref{_gene_torus_no_Hodge_Lemma_} indicates,
the set of tori satisfying \eqref{_very_gene_Equation_}
is also a complement of a countable union of 
positive-codimension subvarieties. We use the
weaker version of generic (\ref{_gene_tori_Definition_})
because it is sufficient for our purposes.


\section{Subvarieties and vector bundles on a generic torus}


\subsection{Subvarieties of generic torus}

\proposition\label{_subva_trivi_Proposition_}
Let $T$ be a generic compact torus, and 
$Z\subset T$ a closed complex subvariety. Then $\dim Z=0$
ot $Z=T$. 

\hfill

{\bf Proof:} Replacing $Z$ with its irreducible
component, we may assume that $Z$ is irreducible. 
Assume that $Z$ has smallest
dimension among all proper subvarieties of $T$.
Then the singular set  $\Sing Z$ is 0-dimensional.
Since $T$ is generic, it 
has no Hodge cycles of dimension 1 and 2.
Therefore $\dim_\C Z >2$

The sheaf $\Omega^1(T)$ is globally generated;
therefore, the same is true for $\Omega^1(Z)$. 
Pick $k$ holomorphic 1-forms $\tau_1, ..., \tau_k$,
such that $\gamma:= \tau_1\wedge ... \wedge \tau_k$
is a non-zero section of the canonical class of $Z$.
Unless the sheaf $\Omega^1(Z)$ is  trivial outside of
 $\Sing Z$, $\gamma$ has non-trivial zero divizor.
This contradicts our assumption that 
$Z$ has minimal dimension among
all subvarieties. Therefore, $\tau_1, ..., \tau_k$
are linearly independent everywhere. These sections
trivialize $\Omega^1(Z)$ outside of $\Sing Z$.

Let $K:= \ker(i^*:\; \Omega^1(T) \arrow \Omega^1(Z))$,
where $i^*$ denotes the pullback map. Since 
$\Omega^1(Z)$ is trivialized, we have $\dim \Omega^1 Z=k$, and
$\dim K = n-k$. Consider $K \subset \Omega^1 T$ as a 
distribution. Then $K$ is integrable, and $Z$ 
is its integral variety. Passing to a covering
of $T$, we find that $Z$ is determined as a 
common zero set of a system of linear equations.
Therefore, $Z$ is flat, and hence it is a subtorus.
This is impossible, because $T$ is generic. \endproof

\subsection{Stable bundles on generic tori}

We follow \cite{_Verbitsky:cohe_}  (see also
\cite{_Voisin_}). 

\hfill

\proposition\label{_gene_stable_Proposition_}
Let $T$ be a generic compact K\"ahler torus, and 
$F$ a stable reflexive\footnote{A sheaf is called
{\bf reflexive} is a natural map $F\arrow \Hom(\Hom(F, \calo_T),\calo_T)$
is isomorphism.} sheaf on $T$. Then $F$ is a line bundle.

\hfill

{\bf Proof:}
In \cite{_Bando_Siu_}, Bando and Siu construct
a canonical Hermitian Yang-Mills  connection $\nabla$
on the non-singular part of any stable
reflexive sheaf on a compact K\"ahler
manifold $X$, with $L^2$-integrable
curvature $\Theta$. Such a connection
(called admissible Yang-Mills) is unique. 
In \cite{_Bando_Siu_} 
it is also proven that the $L^1$-integrable form $\Tr(\Theta\wedge\Theta)$
(considered as a current on the K\"ahler manifold $X$) is closed 
and represents the  cohomology class 
\begin{equation}\label{_curva_Chern_Equation_} 
  \Tr [\Theta\wedge \Theta] = c_2(F) - \frac{r-1}{r} c_1^2(F)
\end{equation}
Consider the standard Hodge operator
$\Lambda:\; \Lambda^{p,q}(T) \arrow \Lambda^{p-1,q-1}(T)$ 
on differential forms. 
Since $\nabla$ is Yang-Mills, we have
\begin{equation}\label{_Lubke_Equation_} 
\Lambda^2 (\Theta^2) \geq 0,
\end{equation}
and the equality is reached only if $\Theta=0$.\footnote{This 
is the celebrated L\"ubke inequality \cite{_Lubke_},
which implies flatness of stable bundles with zero Chern
classes and the Bogomolov-Miyaoka-Yau 
inequality; see  \cite{_Bando_Siu_}.}
Comparing \eqref{_curva_Chern_Equation_} 
and \eqref{_Lubke_Equation_},
and using $c_1(F), c_2(F)=0$, we obtain
that $\Theta=0$, that is, $\nabla$ is flat.
Since $F$ is reflexive, it is 
non-singular in codimension 2 (see \cite{_OSS_}). Therefore,
$\nabla$ has no local monodromy around singularities
of $F$. This implies that $(F,\nabla)$ can be extended to a flat holomorphic
bundle on $T$. Any reflexive sheaf is normal,
that is, equal to a pushforward from any open set
$T\backslash Z$, with $\codim Z> 1$,
(see \cite{_OSS_}). Therefore, $F$ is smooth.

We proved that $F$ is a smooth holomorphic 
bundle admitting a flat unitary connection. Denote by 
$\chi$ the corresponding representation of $\pi_1(T)$.
Since $F$ is stable, $F$ cannot be split
onto a direct sum of subsheaves. Therefore
$\chi$ is irreducible. Since $\pi_1(T)$
is abelian, its irreducible representation
is necessarily 1-dimensional. We proved
that $F$ is a line bundle. \endproof

\subsection{Reflexive sheaves on generic tori}

We follow \cite{_Verbitsky:cohe_}.

\hfill

\proposition\label{_gene_refle_Proposition_}
Let $T$ be a generic compact K\"ahler torus, and 
$E$ a reflexive sheaf on $T$. Then $E$ is a bundle.

\hfill

{\bf Proof:} Since $T$ has no non-trivial
integer $(1,1)$-cycles,
the first Chern class of all sheaves 
on $T$ vanishes. Therefore, all
coherent sheaves on $T$ 
 are semistable.  Consider the Jordan-Holder
filtration
\[
0 = E_0 \subset E_1 \subset ... \subset E_n = E
\]
with all the subfactors $E_i/E_{i-1}$ stable.
The reflexizations $(E_i/E_{i-1})^{**}$ are 
all smooth, by \ref{_gene_stable_Proposition_}. 
Replacing $E_i$ by its reflexization $E_i^{**}\subset E$, we 
may assume that all $E_i$ are reflexive.
Using induction, we may assume also that $E_{n-1}$
is smooth. 

We have an exact sequence
\[
0 \arrow E_{n-1} \arrow E \arrow E_n / E_{n-1} \arrow 0
\]
with $E_{n-1}$ smooth, $E$ reflexive, and $F = E_n / E_{n-1}$
a stable sheaf having (as we have shown above) a smooth
reflexization.

Then $E$ is given by a class $\nu \in \Ext^1(F, E_{n-1})$.
Consider the exact sequence
\[ 0 \arrow F \arrow F^{**} \arrow C \arrow 0
\]
where $C$ is a torsion sheaf (cokernel of the reflexization map).
This gives a long exact sequence
\begin{equation}\label{_long_refle_on_torus_Equation_}
\Ext^1(F^{**}, E_{n-1}) \arrow \Ext^1(F, E_{n-1})
\stackrel \delta \arrow \Ext^2(C, E_{n-1}) 
\end{equation}
The kernel of $\delta$ in \eqref{_long_refle_on_torus_Equation_}
corresponds to all extensions 
$\gamma\in \Ext^1(F, E_{n-1})$
with a reflexization
isomorphic to an extension of $F^{**}$ with $E_{n-1}$.
Clearly, such extensions are reflexive only if $C=0$.
To prove that $C=0$, it suffices to show that $\delta(\nu)=0$.
However, $C$ is a torsion sheaf, and by 
\ref{_subva_trivi_Proposition_}
 its support $Supp(C)$
is a finite set. By Grothendieck's duality
(\cite{_Hartshorne:LC_}, Theorem 6.9),
the group $\Ext^2(C, E_{n-1})$
vanishes, for $\codim \Sup C >2$. Therefore,
$E$ is smooth. This proves 
\ref{_gene_refle_Proposition_}.
\endproof


\section{$SO(2n)/U(n)$-twistor space}


Let $(V, g)$ be a Euclidean vector space, $\dim_\R V = 2n$,
and $S$ the set of all Hermitian structures on $V$ compatible
with $g$. Clearly, $S\cong SO(2n)/SU(n)$. The space $S$ is
a K\"ahler symmetric space, which can be seen from the
following argument.

Consider the space $V_\C:= V\otimes \C$, and let $S'$ 
be the Grassmanian space of all complex subspaces
$V_1 \subset V_\C$ such that $\dim_\C V_1 = n$,
and $g_\C\restrict {V_1}=0$, where $g_\C$ is 
the complexification of $g$. Since $g$ is positive
definite, $V_1 \cap V =0$. Therefore, the natural
projection $V_\C \arrow V$, $v \arrow \Re(v)$ 
induces an isomorphism $V_1 \cong V$. The
space $V_\C$ is equipped with a Hermitian structure
as a complexification of Euclidean space. This gives
a Hermitian structure on $V$. It is easy to check
that this Hermitian structure is compatible with the
Euclidean structure on $V$. Conversely, any such
structure defines a subspace $V^{1,0}\subset V_\C$ 
of vectors of Hodge type (1,0), which is obviously
isotropic. We obtain 
a bijection $S'\cong S$. However, $S'$ is
by construction a K\"ahler symmetric space.

\hfill

\definition\label{_S_Definition_}
We call $S= SO(2n)/U(n)$ {\bf the isotropic Grassman manifold}.

\hfill

Given a torus $T:= V / \Z^{2n}$, we may identify $S$ with the
space of complex structures on $T$ compatible with the metric $g$.
Consider the product $T\times S$, equipped with the complex structure
as follows. 

Let $(t, s)\in T\times S$ be a point. As we have explained, $s$
defines a complex structure $I_s$ on $T$. Let 
${\cal I}:\; T_tT \times T_sS\arrow T_tT \times T_sS$
act as $I_s$ on $T_tT$ and as $I_S$ on $T_sS$, where
$I_S$ is the standard complex structure on $S$.
Clearly, ${\cal I}$ is an almost complex structure
on $T\times S$. 

\hfill

\proposition\label{_twi_inte_Proposition_}
In the above assumptions, ${\cal I}$ is integrable. 

\hfill

{\bf Proof:} Let $B^-$ be the universal bundle on
\[S= \{ V^{1,0} \subset V_\C \ \ | \ \  g\restrict{V_1}=0\},
\]
with the fiber of $B^-$ at a point $V_1\subset V_\C$ identified
with $V^{1,0}$. By construction, $B^-\subset V_\C \otimes \calo_S$ 
is holomorphic. We have an exact sequence
\begin{equation} \label{_exa_sequa_B^+_Equation_}
0 \arrow B^- \arrow V_\C \otimes \calo_S \stackrel {\kappa} \arrow B^+ \arrow 0
\end{equation}
Consider the lattice $L\subset V$, $T= V /L$,
and let $\alpha_1, ..., \alpha_{2n}$ be its generators.
Consider the corresponding sections 
$\alpha_1\otimes 1 , ..., \alpha_{2n}\otimes 1$
of $V_\C \otimes \calo_S$. 
Since $B^-\restrict s$ is an isotropic
subspace of $V_\C= V_\C \otimes \calo_S\restrict s$, 
for all $s \in S$, no non-trivial linear combination of
$\alpha_i\otimes 1$ lies in $B^-$. Therefore, 
the sections $\kappa(\alpha_i\otimes 1)\in B^+$
are linearly independent over $\R$, and
the quotient $B^+/\langle \kappa(\alpha_i\otimes 1)\rangle$
is a compact torus, at every point $s\in S$.
The total space $\Tw(T)$ of this quotient is a
holomorphic fibration over $S$, and its total
space is naturally identified with $(T\times S, {\cal I})$.
This proves \ref{_twi_inte_Proposition_}. \endproof

\hfill

\definition
Let $(T,g)$ be a compact complex torus equipped
with a flat metric. The complex manifold
\[ \Tw(T):= (T\times S, {\cal I})\] is called
{\bf the $SO(2n)/U(n)$-twistor space of $(T,g)$},
or simply twistor space of $T$.
Clearly, $\Tw(T)$ is equipped with a holomorphic
projection $\pi:\; \Tw(T)\arrow S$, and
its fibers are identified with $(T, I_s)$,
where $I_s, s\in S$ are complex structure
operators on $T$ compatible with $g$.


\section[Twistor transform for $SO(2n)/U(n)$-twis\-tor space]{Twistor transform for $SO(2n)/U(n)$-twis\-tor \\ space}


\subsection{Twistor sections for $\Tw(T)$}

\definition
Let $T$ be a compact torus, and \[ \Tw(T) \stackrel \pi\arrow S\]
its $SO(2n)/U(n)$-twistor space. A twistor section
of $\Tw(T)$ is a holomorphic map 
$S\stackrel i \hookrightarrow \Tw(T)$
such that the composition
\[
S\stackrel i \hookrightarrow \Tw(T)\stackrel \pi\arrow S
\]
is identity. We denote the space of twistor sections by $\Sec(T)$.
For any point $t\in T$, the map $I \arrow (I, t)$
gives a holomorphic section of $\pi$. Such sections
are called {\bf horizontal twistor sections}.

\hfill

\proposition\label{_twi_sec_by_two_points_Proposition_}
In the above assumptions, let $\tau:\; S \arrow \Tw(T)$
be a twistor section, and $I, J\in S$ points in 
\[S= \{ V^{0,1} \subset V_\C \ \ | \ \  g\restrict{V_1}=0\},
\]
such that the corresponding $V^{0,1}$-spaces do not intersect:
\begin{equation}\label{_V^0,1_dont_intersect_Equation_}
V^{0,1}_I \cap V^{0,1}_J =0.
\end{equation}
Then, in a neighbourhood $U$ of $[\tau]\in\Sec(T)$,
every section $\tau'$ is uniquely determined by
$\tau'(I)$, $\tau'(J)$. This defines an
open embedding $U\hookrightarrow (T,I)\times (T,J)$.

\hfill

{\bf Proof:} By construction, $\Tw(T)$ is a quotient of 
the total space $\Tot B^+$ by a lattice. Since 
$S$ is simply connected, every section
$\tau\in \Sec(T)$ is lifted to a 
holomorphic section of the projection
\[
\Tot B^+ \stackrel{\tilde \pi} \arrow S.
\]
Therefore, to prove \ref{_twi_sec_by_two_points_Proposition_}
it suffices to show that a section of $B^+$ is uniquely
determined by its values at $I$ and $J$. This is implied
by the following

\hfill

\lemma\label{_section_B^+_determined_Lemma_}
Let $S$ be an isotropic Grassmanian
(\ref{_S_Definition_}),  $B^+$ be the universal
bundle (see \eqref{_exa_sequa_B^+_Equation_}), and 
$I$, $J\in S$ points satisfying 
$V^{0,1}_I \cap V^{0,1}_J=0$. 
Consider the restriction maps
\[
H^0(B^+) \stackrel {r_I} \arrow B^+\restrict I, \ \ 
H^0(B^+) \stackrel {r_J} \arrow B^+\restrict J. 
\]
Then 
\begin{equation} \label{_restict_B^+_iso_Equation_}
H^0(B^+) \stackrel {r_I\oplus r_J} \arrow B^+\restrict I\oplus B^+\restrict J.
\end{equation}
is an isomorphism of vector spaces. In other words, a section
of $B^+$ is uniquely determined by its values at $I$, $J$.

\hfill

{\bf Proof:} Consider the exact sequence
\[
0 \arrow B^- \arrow V_\C \otimes \calo_S \stackrel {\kappa} \arrow B^+ \arrow 0
\]
\eqref{_exa_sequa_B^+_Equation_}. As \ref{_coho_B^-_vanishes_Lemma_}
(below) implies,
$H^0(B^-)=H^1(B^-)=0$. From the corresponding long
exact sequence, we obtain that the induced
map 
\begin{equation}\label{_H^0B^+_is_V_C_Equation_}
H^0(V_\C \otimes \calo_S) \stackrel\kappa \arrow H^0(B^+)
\end{equation}
is an isomorphism. Therefore, $H^0(B^+)$ is naturally
identified with $V_\C$, and $\dim H^0(B^+)= 2 \rk B^+$.
We obtain that the spaces on the left and the 
right hand side of \eqref{_restict_B^+_iso_Equation_} 
have the same dimension.

To prove that \eqref{_restict_B^+_iso_Equation_}
is an isomorphism of vector spaces, it suffices
to show that $\ker r_I \cap \ker r_J=0$. 
This is equivalent to $V^{0,1}_I \cap V^{0,1}_J=0$,
because under the isomorphism \eqref{_H^0B^+_is_V_C_Equation_}
$\ker r_I$ corresponds to $V^{0,1}_I$ and
$\ker r_J$ corresponds to $V^{0,1}_J$.
We reduced \ref{_section_B^+_determined_Lemma_}
and \ref{_twi_sec_by_two_points_Proposition_}
to the following algebro-geometric lemma.

\hfill

\lemma\label{_coho_B^-_vanishes_Lemma_}
Let $S$ be the isotropic Grassmanian
(\ref{_S_Definition_}) and 
\begin{equation}\label{_main_equa_2nd_Equation_}
0 \arrow B^- \arrow V_\C \otimes \calo_S \stackrel {\kappa} \arrow B^+ \arrow 0
\end{equation}
the exact sequence defined above. Then $H^0(B^-)=H^1(B^-)=0$.

\hfill

{\bf Proof:} We use the Kodaira-Nakano vanishing theorem.
Since $B^-$ is a holomorphic sub-bundle of a trivial
bundle, it is negative: $\Theta_{B^-}\leq 0$,
where $\Theta_{B^-}\in \Lambda^{1,1}(S) \otimes \End(B^-)$
is the curvature of $B^-$ (\cite{_Griffi_Harri_}).
The form $\Theta_{B^-}$ can be computed explicitly
as follows:
\[ \Theta_{B^-}= A \wedge A^\bot,
\]
where 
\begin{equation}\label{_A_where_Equation_}
A \in \Lambda^{1,0}(S) \otimes \Hom(B^-, B^+)
\end{equation}
is the second form of the exact sequence
\eqref{_main_equa_2nd_Equation_}, and 
\[ A^\bot \in \Lambda^{0,1}(S) \otimes \Hom(B^+, B^-)\]
its Hermitian conjugate (see \cite{_Griffi_Harri_}
for details). The second form $A$ can be computed
explicitly, as follows. 

The tangent bundle $T(S)$ is naturally identified
with $\Hom(B^-, B^+)$, because $S$ is a Grassmanian,
and a tangent vector to a Grassmanian of planes
$V_0\subset V_1$ is identified with $\Hom(V_0, V_1/V_0)$.
Under the identification 
\begin{equation}\label{_Lambda^1(S)_Equation_}
\Lambda^{1,0}(S)\cong \Hom(B^+, B^-),
\end{equation}
$A$ becomes an identity operator.
\begin{multline*}
A \in \Lambda^{1,0}(S) \otimes \Hom(B^-, B^+) \\=
  \Hom(B^+, B^-)\otimes \Hom(B^-, B^+) \\ = \End(\Hom(B^-, B^+)).
\end{multline*}
This implies that $\Theta_{B^-}= A \wedge A^\bot$
is strictly negative. Now, by Kodaira-Nakano theorem,
$H^i(B^-)=0$ for all $i< \dim S$. This proves 
\ref{_coho_B^-_vanishes_Lemma_}. We finished the
proof of \ref{_twi_sec_by_two_points_Proposition_}.
\endproof

\hfill

\remark \label{_complexi_Remark_}
Consider an anticomplex involution 
$\iota:\; S\arrow S$, $I \arrow -I$.
Then $\iota \times Id$ is an anticomplex involution
of $\Tw(T)$. This involution maps twistor sections to twistor
sections. Abusing notation, we 
denote the corresponding involution of
$\Sec(T)$ by $\iota$.

There is a natural real analytic embedding
$T\stackrel \gamma \hookrightarrow \Sec(T)$, $t\arrow \{t\}\times S$ mapping
$t$ to the corresponding holomorphic section.
Clearly, $\gamma(T)$ is a fixed point set of
the anticomplex involution $\iota:\; \Sec(T)\arrow \Sec(T)$.
By \ref{_twi_sec_by_two_points_Proposition_},
$\Sec(T)$ is locally in a neighbourhood
of $\gamma(T)$ identified with a
neighbourhood of diagonal in 
$(T,I)\times (T, -I)$.
We shall think of 
$\Sec(T)$ as of complexification 
of its totally real submanifold $T=\gamma(T)$.

\subsection{Twistor transform for $\Tw(T)$}

Let $(B, \nabla)$ be a flat bundle on $T$, and
$(\sigma^* B, \sigma^*\nabla)$ its pullback to
$\Tw(T)$. Clearly, $(\sigma^* B, \sigma^*\nabla)$
is flat, and therefore holomorphic. This defines
a functor from the category of flat bundles on $T$
to the category of holomorphic bundles on $\Tw(T)$:
\begin{equation}\label{_dire_twi_fu_Equation_}
\Fl(T) \arrow \HolBun(\Tw(T)).
\end{equation}
We call this functor {\bf the direct twistor transform}.
It turns out that \eqref{_dire_twi_fu_Equation_}
is invertible. 

The following theorem is adapted from \cite{_NHYM_},
where a similar result was proven for a twistor space of 
a hyperk\"ahler manifold. 

\hfill

\theorem \label{_iva_twi_Theorem_}
Let $T$ be a torus, $\Tw(T)\stackrel \pi \arrow S$ its
 $SO(2n)/U(n)$-twistor space, and
$B$ a holomorphic vector bundle on $\Tw(T)$. Assume that
\begin{equation}\label{_trivial_on_horiz_Equation_}
B\restrict {S\times \{t\}\subset S\times T =\Tw(T)}
\text{\ \it \scriptsize is trivial for all\ }t\in T, 
\end{equation}
that is, $B$ is trivial
on all horizontal sections of $\pi$.
Then $B$ is obtained as a twistor transform of
a flat bundle $(B_\fl, \nabla)$ on $T$.
Moreover, $(B_\fl, \nabla)$ is unique, and 
determines an equivalence 
\[
\HolBun_0(\Tw(T))\arrow \Fl(T)
\]
of the category $\HolBun_0(\Tw(T))$
of holomorphic bundles on $\Tw(T)$ satisfying 
\eqref{_trivial_on_horiz_Equation_} and the category
of flat bundles on $T$.

\hfill

We prove \ref{_iva_twi_Theorem_} in Subsection 
\ref{_real_analy_holo_on_S_Subsection_}.

\subsection{Real analytic differential operators and
holomorphic vector bundles on twistor spaces}
\label{_real_analy_holo_on_S_Subsection_}

We work in assumptions of \ref{_iva_twi_Theorem_}.
Let $B\in \HolBun_0(\Tw(T))$ be a holomorphic vector
bundle which is trivial on all horisontal sections
of twistor projection. Denote by $U_B\subset \Sec(T)$
the space of all twistor sections 
$S\stackrel \phi\hookrightarrow \Tw(T)$ for which 
$B\restrict{\phi(S)}$ is trivial. 

\hfill

\lemma\label{_U_b_open_Lemma_}
In the above assumptions, $U_B$ is open in $\Sec(T)$.

\hfill

{\bf Proof:} We need to show that a small deformation of a trivial
bundle on $S$ is again trivial. It is well known that
deformations of vector bundles are classified by
$\Ext^1(B,B)= H^1(\End(B))$ (\cite{_Kobayashi_}). If $B$ is a trivial
bundle, $\End(B)$ is also trivial. Therefore, 
\ref{_U_b_open_Lemma_} is implied by 
\begin{equation}\label{H^1(calo)_Equation_}
H^1(\calo_S)=0.
\end{equation}
{}From \eqref{_Lambda^1(S)_Equation_} 
and ampleness of $B^+$, we obtain that 
$TS$ is ample (this is also clear because $S$ 
is symmetric). Therefore, the canonical class $K$ of $S$ 
is negative. By Serre's duality, 
\[ H^1(\calo_S) = H^{\dim S-1}(K)^*.
\]
By Kodaira-Nakano theorem, $H^i(K)=0$ for all $i< \dim S$,
because $K$ is negative. This proves \ref{_U_b_open_Lemma_}.
\endproof

\hfill

Consider the evaluation map $S\times U_B\stackrel {ev}\arrow \Tw(T)$,
and let $ev^* B$ be the pullback of $B$. Then $ev^*B$ is trivial on
all $S\times \{s\}\subset S\times U_B$, and therefore, the
pushforward $B_\C= \pi_*(ev^* B)$ is a holomorphic
vector bundle on $U_B$.

One should think of $U_B$ as of complexification $T_\C$ of
$T$ (\ref{_complexi_Remark_}). Then $B_\C$ becomes
a complexification of a real analytic bundle on 
$T$ underlying $B$. The holomorphic structure operator 
on $B$ can be interpreted as a holomorphic
differential operator on $B_\C$, as follows.

Fix $I \in S$, and let $T_I \subset T U_B$ 
be a tangent hyperplane field to a holomorphic
foliation $\c L$ on $\Sec(T)$ with leaves 
\[ L_t:= \{ \phi:\; S \arrow \Tw(T) \ \ |\ \ \phi(I) = t \}\]
parametrized by $t\in T$. We can write
$L_t = ev^{-1}_I(t)$. Clearly, 
$B_\C\cong ev^*_I(B\restrict{(T,I)})$.
Therefore, $B_\C$ is trivialized along the leaves
of $\c L$. This defines a trivial 
holomorphic conection along
the leaves of $\c L$:
\begin{equation}\label{_holo_str_via_leaves_Equation_}
\nabla_I:\; B_\C \arrow B_\C \otimes T^*_I.
\end{equation}
Restricting \eqref{_holo_str_via_leaves_Equation_} to
$\gamma(T)\subset \Sec(T)$ (see \ref{_complexi_Remark_}),
we obtain that $T^*_I$ becomes $\Lambda^{0,1}(T)$,
and \eqref{_holo_str_via_leaves_Equation_}
becomes the holomorphic structure operator
\[ \bar\6:\; B\arrow B\otimes \Lambda^{0,1}(T).\]

Writing this operator as in
\eqref{_holo_str_via_leaves_Equation_}
allows us to study the dependence of $\bar\6$ from $I$ in
a more explicit way. 

Since $U_B$ is an open neighbourhood of its 
totally real submanifold $\gamma(T)$, 
it contains a Stein neighbourhood of
$\gamma(T)$ as Grauert's theorem implies.
Replacing $U_B$ by a smaller neighbourhood
if necessary, we may assume that $U_B$ is Stein.
Then $B_\C$, $B_\C \otimes T^*_I$ are globally
generated, and the map \eqref{_holo_str_via_leaves_Equation_}
can be interpreted as a map of $\calo_{U_B}$-modules:
\[
\nabla_I:\; H^0(B_\C) \arrow H^0(B_\C\otimes T^*_I).
\]
Clearly, \eqref{_holo_str_via_leaves_Equation_}
depends holomorphically from the parameter $I\in S$. 
Therefore, \eqref{_holo_str_via_leaves_Equation_}
can be interpreted as a map
\begin{equation}\label{_holo_stru_from_I_Equation_}
\nabla_I:\; H^0(B_\C)\otimes_\C \calo_S \arrow H^0(B_\C)\otimes_\C B^+.
\end{equation}
(we use the identification of $T^*_I=\Lambda^{0,1}(T)$ with $B^+$
obtained earlier). 

Let $\Gamma$ be a vector space, possibly 
infinite-dimensional.
Consider the exact sequence 
\begin{equation}\label{_exa_tensored_by_Gamma_Equation_}
0 \arrow B^- \otimes \Gamma \arrow V_\C \otimes \calo_S\otimes \Gamma\stackrel {\kappa} \arrow B^+\otimes \Gamma \arrow 0
\end{equation}
(this is \eqref{_exa_sequa_B^+_Equation_} tensored by
$\Gamma$). Since $H^0(B^-) = H^1(B^-)=0$,
the long exact sequence associated with
\eqref{_exa_tensored_by_Gamma_Equation_}
gives 
\[ \Hom (\Gamma_1\otimes \calo_S, \Gamma_2\otimes B^+)\cong
   \Hom (\Gamma_1\otimes \calo_S, \Gamma_2\otimes V_\C \otimes \calo_S),
\]
for any vector spaces $\Gamma_1$, $\Gamma_2$.
Therefore, the map \eqref{_holo_stru_from_I_Equation_}
lifts, uniquely, to a holomorphic map
\begin{equation}\label{_nabla_S_defini_Equation_}
\nabla_S:\; H^0(B_\C) \otimes \calo_S 
  \arrow H^0(B_\C)\otimes V_\C \otimes \calo_S 
  = H^0(B_\C\otimes_{\calo_{U_B}} \Omega^1(U_B))
\end{equation}
It is easy to check that $\nabla_S$ is in fact a 
holomorphic connection on $B_\C$. Restricting
$\nabla_S$ to $\gamma(T)\subset U_B$, we
obtain a real analytic connection operator
$\nabla$ on $B$. By construction,
the $(0,1)$-part of $\nabla$
is the holomorphic structure operator on $B$.

We are going to show that $\nabla$ is flat.
Indeed, the $(0,2)$-part of the curvature $\Theta$ of
$\nabla$ vanishes for all $I\in S$, because $\bar\6=0$.
If we replace $I$ by $-I$, the $(2,0)$-forms become
$(0,2)$-forms and vice versa. Therefore, 
the $(2,0)$-part of the curvature $\Theta$
also vanishes. We obtain that 
$\Theta$ is of type $(1,1)$ with
respect to all $I\in S$. In other words,
$\Theta$ is invariant under the natural
$U(1)$-action induced by the complex structure.
It is easy to see that such $U(1)$ generate
$SO(2n)$. Therefore, $\Theta$ is an $SO(2n)$-invariant
2-form on 2n-dimensional space. Using the standard
invariants theory, we infer that $\Theta=0$.
Therefore, $\nabla$ is flat. From the
construction it is clear that $\nabla$ is unique.

We obtained a functor from $\HolBun_0(T)$
to the category of flat bundles, which is
obviously inverse to the twistor transform.
This proves \ref{_iva_twi_Theorem_}.


\section{Flat connections on semistable bundles}


\subsection{Deformation theory and flat bundles}

\definition
Let $T$ be a compact complex torus. Given a holomorphic bundle
$B$ on $T$, we say that $B$ is {\bf polystable} if $B$
is a direct sum of line bundles.

\hfill

\definition
Let $T$ be a compact complex torus with flat K\"ahler
metric, and $\Tw(T)$ the corresponding twistor space.
Fix a point $I\in S$. Denote by $\Fl_I$ the category
of flat bundles $B$ on $\Tw(I)$ such that $B\restrict{(T,I)}$
is polystable, and every irreducible subquotient of
$B$ admits a Hermitian metric. 

\hfill

We are interested in the category $\Fl_I$ because of the
following theorem.

\hfill

\theorem\label{_restri_from_Fl_I_equi_Theorem_}
Let $T$ be a compact torus, $\Tw(T) \stackrel \pi \arrow S$
its twistor space, and $I, J$ points which
satisfy $V_I^{0,1}\cap V^{0,1}_J=0$
\eqref{_V^0,1_dont_intersect_Equation_}.
Consider the restriction functor from 
$\Fl_I$ to the category $\Bun(T,J)$
of holomorphic bundles on $(T, J) \subset \Tw(T)$:
\begin{equation}\label{_restri_from_Fl_I_Equation_}
\Fl_I \stackrel \Phi\arrow \Bun(T,J).
\end{equation}
Then \eqref{_restri_from_Fl_I_Equation_}
is equivalence of categories.

\hfill

We prove \ref{_restri_from_Fl_I_equi_Theorem_}
in Subsection \ref{_abe_cate_Subsection_}. In the present Section,
we prove a weaker form of \ref{_restri_from_Fl_I_equi_Theorem_}.
A similar result is proven in \cite{_Verbitsky:cohe_}.

\hfill

\proposition \label{_B_restriction_Fl_I_Proposition_}
In assumptions of \ref{_restri_from_Fl_I_equi_Theorem_},
let $B$ be a holomorphic bundle on $(T, J)$. 
Then $B=\Phi(\tilde B)$, for some flat bundle
$\tilde B\in \Fl_I$. Moreover, $\Phi$
induces a bijection on the set of 
isomorphism classes of the objects
of corresponding categories.

\hfill

{\bf Proof:}
We use the following lemma

\hfill

\lemma\label{_exists_nu_flat_Lemma_}
Let $T$ be a compact complex torus,
$B_1, ... B_n$ flat holomorphic Hermitian vector bundles
and $B$ a holomorphic vector bundle with a filtration
\[ 0 = E_0\subset E_1 \subset ...\subset E_n= B,\]
such that $E_i / E_{i-1} \cong B_i$. 
Then the holomorphic
structure on $B$ can be obtained as follows. 
Identify $B$ with $B_{gr} := \oplus B_i$
as $C^\infty$-bundle. Then there is a 
cohomology class $\nu$
\begin{equation}
\nu \in \bigoplus_{i>j} \Ext^1(B_i, B_j) \subset \Ext^1(B_{gr}, B_{gr})
\end{equation}
such that the holomorphic structure operator in $B$
is written as
\begin{equation}\label{_holo_via_nu_Equation_}
\bar \6 = \bar\6_{gr} + \nu_0, 
\end{equation}
where
$\bar\6_{gr}$ is the holomorphic structure
operator on $B_{gr}$, and \[ \nu_0 \in \Lambda^{0,1}(\End(B_{gr}))\]
denotes the harmonic (and hence, flat) representative of $\nu$. 

\hfill

{\bf Proof:}
Write the holomorphic structure operator on $B$ as
$\bar\6 = \bar\6_{gr} +\tilde \nu$,
where $\tilde\nu$ is a (0,1)-form with values
in \[ \oplus_{i>j} \Lambda^{0,1}(T, \Hom(B_i, B_j)).\]
The form $\tilde\nu$ satisfies the Maurer-Cartan
equation
\begin{equation}\label{_MC_Equation_} 
  \bar\6^2 = \bar\6_{gr}(\tilde\nu) + 2 \tilde\nu\wedge\tilde\nu =0
\end{equation}
 Every automorphism $g\in \End B_{gr}$
acts on $\tilde \nu$ as 
$\tilde\nu\arrow g(\tilde\nu) + \bar\6_{gr}(g)$
(this is the well-known gauge action).
To produce $\nu$ with the properties described
in \ref{_exists_nu_flat_Lemma_}, we
need to find a correct gauge transform.

Consider the group $(\C^*)^n$ acting on
$B_{gr}$ by diagonal automorphisms, in such a way that
the $i$-th component $\alpha_i$ of $(\C^*)^n$ acts trivially
on $B_j\subset B_{gr}$ for $i\neq j$, and
as a multiplication by $\alpha_i$ on $B_i$.

We shall write the action of $(\C^*)^n$
on $\Lambda^{0,1}(\End(B_{gr})$ as follows.
Let 
\begin{align*} \tilde \nu &\in \Lambda^{0,1}(\End(B_{gr}) = 
   \oplus_{i,j}\Lambda^{0,1}(B_i, B_j),  \\ 
   \tilde\nu &:=  \sum_{i, j}\tilde\nu_{ij}, \ \ \ 
   \tilde\nu_{ij}\in \Lambda^{0,1}(B_i, B_j)
\end{align*}
If $\alpha \in(\C^*)^n$, $\alpha = \prod_i \alpha_i$,
then 
\[ \alpha(\tilde \nu) = \sum_{i, j} \alpha_i \alpha_j^{-1}
   \tilde\nu_{ij}
\]
The group $(\C^*)^n$ acts in this fashion on the solutions of
Maurer-Cartan equation, and maps every solution to an equivalent
one. If $\alpha_j \gg \alpha_i$ for all $i>j$, then
$\alpha$ maps 
\[ \tilde\nu\in\oplus_{i>j} \Lambda^{0,1}(T, \Hom(B_i, B_j))
\]
to a form which is arbitrarily small.

Consider the local deformation space $\Def(B_{gr})$ for
$B_{gr}$, constructed in \cite{_Siu_Trautmann_}.
The above argument implies that every
neighbourhood of the point $[B_{gr}]\in \Def(B_{gr})$
contains a bundle which is isomorphic to $B$.

Let $E$ be a holomorphic vector bundle over a compact
K\"ahler manifold. The local deformation space $\Def(E)$
can be constructed explicitly in terms of Massey products 
as follows. 

One can define the Massey products as obstructions
to constructing a solution of the Maurer-Cartan equation
(see e.g. \cite{_Babenko_Taimanov_}, or
\cite{_May_}, \cite{_Retakh_} for a more 
classical approach). Locally, $\Def(E)$ is 
embedded to the vector space $\Ext^1(E, E)$,
and the image of this embedding is a germ of all
vectors $\theta\in \Ext^1(E, E)$, such that
$\theta\wedge\theta=0$ and all the higher Massey
products of $\theta$ with itself vanish.

Fix such a vector $\theta\in \Ext^1(E, E)$.
We construct the corresponding 
vector bundle $E_\theta\in \Def(E)$
using the Hodge theory as follows (see e.g. 
\cite{_Verbitsky:Hyperholo_bundles_}).

Let $\theta_0 \in \c H^{0,1}(\Hom(E,E))$
be the harmonic representative of $\theta$. Using
induction, we define
\[ \theta_n := - \frac{1}{2} G_{\bar\6}\sum_{i+j=n-1} \theta_i \wedge \theta_j
\]
where $G_{\bar\6}$ is the Green operator inverting the
holomorphic structure operator
\[ \bar\6:\; \Lambda^{0,k}\otimes E\arrow  
    \Lambda^{0,k+1}\otimes E.
\]
on its image.
The Green operator $G_{\bar\6}$ is compact.
This can be used to show that 
for $\theta$ sufficiently small, the
series $\tilde \theta:= \sum \theta_i$
converges. The vanishing of Massey products
is equivalent to the following condition
\[ \bar\6 \theta_n = 
    - \frac{1}{2} \sum_{i+j=n-1} \theta_i \wedge \theta_j,
\]
which is apparent from the definition given in
\cite{_Babenko_Taimanov_}. In this case, we have
\[ \bar\6 \tilde \theta = 
   - \frac{1}{2} \tilde\theta\wedge\tilde\theta
\]
and $\tilde\theta$ is a solution of the Maurer-Cartan equation
\eqref{_MC_Equation_}.
Therefore, the operator $\bar\6_\theta= \bar\6+\tilde\theta$
satisfies $\bar\6_\theta^2=0$, and by Newlander-Nirenberg
theorem (\cite{_Kobayashi_}, Proposition 4.17, Chapter 1)
this operator defines a holomorphic structure
on $E$. On the deformation 
space $\Def(E) \subset \Ext^1(E, E)$,
the point $\theta$ corresponds to a bundle
$(E, \bar\6_\theta)$.

Now we return to holomorphic bundles over a compact torus
and the proof of \ref{_exists_nu_flat_Lemma_}. We obtain that
$B$ is given by some $\nu\in \Ext^1(B_{gr}, B_{gr})$. Since
$B$ and $\oplus B_{i}$ have compatible filtrations,
by functoriality we may assume that
\begin{equation}\label{_nu_nilpo_Equation_}
\nu \in \oplus_{i>j} \Ext^1(B_i, B_j).
\end{equation}
The higher Massey operations in $\Ext^*(B_{gr}, B_{gr})$
vanish, because the bundle $B_{gr}$ is flat
(the same proof works as was used in \cite{_DGMS_};
see also \cite{_Goldman_Millson_}). Therefore,
$B$ can be reconstructed from $\nu$ 
for any $\nu$ such that the cohomology class 
$\nu\wedge\nu$ vanishes. Pick a harmonic representative
$\nu_0$ of $\nu$. Since $B_{gr}$ is flat, $\nu_0$ is
parallel. Therefore $\nu_0\wedge\nu_0$ is also parallel,
hence harmonic. We obtain that the cohomology
class of $\nu_0\wedge\nu_0$ vanishes if and only
if this form vanishes identically.

Starting from a bundle $B$ with a filtration
satisfying the assumptions of \ref{_exists_nu_flat_Lemma_},
we have constructed a cohomology class 
$\nu \in \oplus_{i>j} \Ext^1(B_i, B_j)$, with
$\nu\wedge\nu=0$. The harmonic representative $\nu_0$ of $\nu$
satisfies $\nu_0\wedge\nu_0=0$. Therefore, $\nu_0$ 
is a solution of Maurer-Cartan equation, and 
$\bar\6_{gr} +\nu_0$ is equivalent to 
the holomorphic structure operator of $B$.
This proves \ref{_exists_nu_flat_Lemma_}.
\endproof

\hfill

Let $\nabla_{gr}$ be the standard 
Hermitian flat connection on $B_{gr}$.
The bundle $(B_{gr}, \nabla_{gr} + \nu_0)$ is also flat,
because $\nu_0$ is by construction parallel.
We obtained the following corollary.

\hfill

\corollary\label{_bun_admits_flat_Corollary_}
Let $B$ be a holomorphic vector bundle on a compact
generic torus $T$, $\dim_\C T>2$.
Then $B$ admits a flat connection 
compatible with the holomorphic structure.

\endproof

\hfill

Return to the proof of 
\ref{_B_restriction_Fl_I_Proposition_}.
We use the notation of \ref{_exists_nu_flat_Lemma_}.
Let $B$ be a holomorphic bundle on $(T,J)$,
and $B_{gr}$ the associated graded bundle of
Jordan-Holder filtration. Then $B_{gr}$ 
is polystable, and hence admits 
a flat Hermitian connection 
(\ref{_gene_stable_Proposition_}). 
Consider the bundle $\tilde B_{gr}$ on $\Tw(T)$, 
obtained from the flat bundle 
$B_{gr}$ via twistor transform, and let
$R^1\pi_*(\End(\tilde B_{gr}))$ be the corresponding  
direct image sheaf on $S$. Clearly,
$R^1\pi_*(\End(\tilde B_{gr}))$ is isomorphic
to a sum of several copies of the vertical
tangent bundle $T_{vert}\Tw(T)\cong B^+$. 
Given a holomorphic section 
$\tilde \gamma \in R^1\pi_*(\End(\tilde B_{gr}))$,
$\gamma^2=0$, we take the fiberwise 
harmonic (hence, flat) representative $\tilde \gamma_0$ 
of $\tilde \gamma$. The  bundle 
$(\tilde B_{gr}, \bar\6_{gr} + \tilde \gamma_0)$
is holomorphic. Consider the 1-form 
$\nu_0\in H^1((T, J),\End(\tilde B_{gr}))$
associated with $B$ as in  \ref{_exists_nu_flat_Lemma_}.

Let $\tilde \nu_0$ be a section of
$R^1\pi_*(\End(\tilde B_{gr}))$ which vanishes
in $I$ and is equal to $\nu_0$ at $J$.
Such a section exists by \ref{_section_B^+_determined_Lemma_}
(we surmise that $R^1\pi_*(\End(\tilde B_{gr}))$
is isomorphic to a sum of several copies of
$B^+$). Clearly, $\tilde B :=(\tilde B_{gr}, \bar\6_{gr} + \tilde \nu_0)$
is a holomorphic bundle with $\tilde B \restrict{(T,I)}=B_{gr}$
and $\tilde B \restrict{(T,J)}=B$. Therefore, $\tilde B$
belongs to $\Fl_I$, and $\Phi(\tilde B)=B$. This bundle
is uniquely determined by $B$, because
$\tilde \nu_0$ is uniquely determined by 
the conditions $\tilde\nu_0\restrict I=0$, 
$\tilde\nu_0\restrict J=\nu$
(\ref{_exists_nu_flat_Lemma_}).  This proves 
\ref{_B_restriction_Fl_I_Proposition_}. \endproof

\hfill

\subsection{Abelian categories of finite length}
\label{_abe_cate_Subsection_}

\definition
Let $\c C$ be an abelian category.
An object $B\in \c C$ is called {\bf simple}
if it has no proper sub-objects:
for all $B'\subset B$, either $B'=B$
or $B'=0$. An object $B$ is called {\bf semisimple}
if $B$ is a direct sum of simple objects.

An abelian category $\c C$ is called {\bf of finite length } if 
every object $B\in \c C$ is a finite extension
of simple objects. Equivalently, $\c C$ 
is of finite length if any increasing or decreasing chain of 
sub-objects of $B$ stabilizes, for all 
$B\in \c C$.\footnote{One also says {\bf $\c C$ 
satisfies the ascending and descending chain condition.}}

A {\bf length} of $B\in \c C$ is a minimal length
of such a chain.

\hfill

\lemma \label{_embe_on_Ext_1_hence_full_faith_Lemma_}
Let $\c C \stackrel \gamma \arrow \c C'$ be a functor of
abelian categories of finite length. Assume that
$\gamma$ induces equivalence on the respective
subcategories of semisimple objects. Assume,
moreover, that $\gamma$ induces a bijection
on the sets of equivalence classes of the objects of
$\c C$, $\c C'$. Then $\gamma$ is equivalence.

\hfill

{\bf Proof:} The proof is obtained via the trivial diagram-chasing
argument. We use an induction by the length of $B$. 
 Let $B, B'\in \c C$ be objects of 
length $\leq n$. Consider an exact sequence
\[ 0 \arrow B_0 \arrow B \arrow B_1 \arrow 0
\]
with $l(B_1)= n-1$ and $B_0$ semisimple. There is a
long exact sequence
\[
0 \to \Hom(B', B_0) \to \Hom(B', B) \to 
\Hom(B', B_1) \to \Ext^1(B', B_0) \to ...
\]
Applying $\gamma$, we obtain a commutative diagram with
exact rows

{\tiny
\minCDarrowwidth0.9pc
\begin{equation}\label{_long_exa_Phi_Ext^1_Equation_}
\begin{CD}
0 @>>> \Hom(B', B_0) @>>> \Hom(B', B) @>>> 
\Hom(B', B_1)@>>> \Ext^1(B', B_0)\\
&& @V{\gamma} VV @V{\gamma}VV  @V{\gamma}VV @V{\gamma}VV \\
0 @>>> \Hom(\gamma(B'), \gamma(B_0)) @>>> \Hom(\gamma(B'), \gamma(B)) @>>> 
\Hom(\gamma(B'), \gamma(B_1)) @>>>\Ext^1(\gamma(B'), \gamma(B_0))
\end{CD}
\end{equation}
}
The maps of $\Ext^1$-groups are isomorphisms
because $\Ext^1$ classify the classes of extensions,
and we know that $\gamma$ induces bijection
on isomorphism classes of objects.
Using induction by length of $B$, $B'$, 
we may assume that all vertical arrows of
\eqref{_long_exa_Phi_Ext^1_Equation_} are isomorphisms,
except, possibly, the map
\[ \Hom(B', B)\stackrel \gamma\arrow \Hom(\gamma(B'), \gamma(B)).\]
Using the snake lemma, we obtain that this map is also
an isomorphism. Therefore, $\gamma$ is equivalence.
\endproof

\hfill

We apply \ref{_embe_on_Ext_1_hence_full_faith_Lemma_}
to prove \ref{_restri_from_Fl_I_equi_Theorem_}.
Simple objects in $\Bun(T,J)$ are line bundles,
as \ref{_gene_stable_Proposition_} implies.
Indeed, simple objects must be stable,
and stable bundles are line bundles.
Clearly, line bunles on $T$ admit flat
Hermitian connections.

As a definition of $\Fl_I$ implies,
simple objects here are line
bundles obtained from flat Hermitian
line bundles $(T,J)$. This implies
that \[ \Phi:\; \Fl_I \arrow \Bun(T,J)\]
is an equivalence on the subcategory of 
semisimple objects. By \ref{_B_restriction_Fl_I_Proposition_},
$\Phi$ induces bijection on equivalence classes of
objects. Now, \ref{_embe_on_Ext_1_hence_full_faith_Lemma_}
implies that it is an equivalence. This proves
\ref{_restri_from_Fl_I_equi_Theorem_}.


\section{Singularities of coherent sheaves on generic torus}


\subsection{Local geometry on $\Tw(T)$}
\label{_local_trivi_Subsection_}

Let $T$ be a compact torus, $\dim T>1$.
Fix a point $x\in T$, and the correspoinding
point in twistor space $(I, x)\in \Tw(T)$.
We shall study the twistor sections passing through
$(I, x)\in \Tw(T)$ in a neighbourhood of 
the horisontal section $s_{x}= S\times\{x\}$.

Let $L, L'\in S$ be complex structures satisfying
$V^{0,1}_L\cap V^{0,1}_I = V^{0,1}_{L'}\cap V^{0,1}_I =0$
(see \eqref{_V^0,1_dont_intersect_Equation_} for details).
Consider a point $(t, L')$ close to $(x, L')$.
By \ref{_twi_sec_by_two_points_Proposition_},
there exists a unique twistor section
$s:\; S \arrow \Tw(T)$ 
passing through $(t, L')$ and $(x, I)$
and close to $s_{x}= \{x\} \times S$.
This section is unique in a small neighbourhood of
$\{x\}\times S$.
We define $\tilde \Psi_I(t):= s(L) \in (T, L)$. 
This map is holomorphic
and invertible in a sufficiently 
small neighbourhood of $(x, L')$.
We think of this map as of isomorphism
\begin{equation}\label{_Psi_on_germs_Equation_}
\Psi_I(L, L'):\; \calo_{x, L} \arrow \calo_{x, L'}
\end{equation}
of the corresponding local rings. Clearly,
$\Psi_I(L, L')$ holomorphically depends on 
$I, L, L'$.

Let now $B$ be a bundle on $\Tw(T)$ obtained from the twistor
transform (\ref{_iva_twi_Theorem_}). This means that
the restriction of $B$ to any horizontal twistor
section $s_x$ is trivial. As we have shown, a small deformation
of a trivial bundle on $S$ is again trivial.
Therefore, the restriction of $B$ to a twistor
section close to $s_x$ is trivial as well. 

This allows one to extend the map
$\Psi_I(L, L'):\; \calo_{x, L} \arrow \calo_{x, L'}$
to the space of germs of holomorphic sections of $B$. 
We obtain an isomorphism
\begin{equation}\label{_Psi_on_germs_of_E_Equation_}
\Psi_I(L, L', B):\; B_{x, L}
  \arrow B_{x, L'}, 
\end{equation}
where $B_{x, L}$ and $B_{x, L'}$ are spaces
of germs of $B\restrict{(T,L)}$ and $B\restrict{(T,L')}$ in $x$.

Denote the infinitesimal neighbourhood
of $s_{x}\backslash {(I, x)}$ in $\Tw(T)\backslash (T, I)$
by 
\[ 
   \Tw(T)_{x, I} = \Spec(\calo_{\Tw(T)\backslash (T, I),s_{x}}).
\]
This space is fibered over $S\backslash I$
with fibers isomorphic to a germ of a smooth
complex manifold.

The maps \eqref{_Psi_on_germs_Equation_}
produce a canonical 
trivialization of $\Tw(T)_{x, I}$
over $S\backslash I$. Denote the trivialization
map by 
\begin{equation}\label{_Phi_triviali_Equation_}
  \Phi:\; \Tw(T)_{x, I} \arrow \Spec \calo_0(\C^n)\times (S\backslash I)
\end{equation}
where $\Spec \calo_0(\C^n)$ is the germ of $\C^n$ in $0$.

Let $\xi:\;  \Tw(T)_{x, I} \arrow \Spec \calo_0(\C^n)$
be the composition of $\Phi$ and the projection
$\Spec \calo_0(\C^n)\times (S\backslash I)\arrow \Spec \calo_0(\C^n)$.
Using the maps $\Psi_I(L, L', B)$
of \eqref{_Psi_on_germs_of_E_Equation_}, we obtain
a trivialization of $B\restrict {\Tw(T)_{x_0, I}}$
over $S \backslash I$. More precisely,
we obtain a bundle $B_S$ over $\Spec \calo_0(\C^n)$, and
a natural isomorphism
\begin{equation}\label{_E_loca_triv_Equation_}
B\restrict {\Tw(T)_{x, I}}\cong \xi^*B_S.
\end{equation}

\subsection{The category $\c C_I$}

\definition\label{_C(I)_Definition_}
Let $T$ be a compact torus, $\Tw(T)$ its 
$SO(2n)/U(n)$-twistor space
The category ${\c C}_I$
is defined as follows. 

An object of ${\c C}_I$
is a coherent sheaf $F$ on $\Tw(T)$
satisfying the following conditions.

\begin{description}
\item[(i)] The reflexization $B:= F^{**}$
belongs to $\Fl_I$.

\item[(ii)] The sheaf $F$ is non-singular outside
of $A\times S\subset \Tw(T)$,
where $A\subset T$ is a finite set.

\item[(iii)] Let $Z:= \{A\} \times  I \subset (T, I) \subset \Tw(T)$
be the finite subset of $(T, I) \subset \Tw(T)$
corresponding to $A$, and 
$j:\; \Tw(T)\backslash Z \hookrightarrow \Tw(T)$
the natural embedding. Then the canonical
homomorphism $F\arrow j_* j^* F$ is an
isomorphism.

\item[(iv)] Let $x\in A$ be a point.
Consider the trivialization
\begin{equation}\label{B_triv_in_defi_Equation_}
B\restrict {\Tw(T)_{x, I}}\stackrel\Psi\arrow \xi^*B_S
\end{equation}
constructed in \eqref{_E_loca_triv_Equation_},
where $B=F^{**}$ is the reflexization of
$F$. Then there is a sheaf $F_S$ on 
$\Spec \calo_0(\C^n)$, equipped with an isomorphism
$F_S^{**}\cong B_S$, and a trivialization
\[
 F\restrict {\Tw(T)_{x, I}}\stackrel\Psi\arrow \xi^*F_S
\]
which is compatible with \eqref{B_triv_in_defi_Equation_}.
\end{description}
The morphisms of ${\c C}_I$
are morphisms of coherent sheaves.

\hfill

\theorem\label{_C(I)_equiv_Theorem_}
Let $T$ be a compact torus,
$\Tw(T)\stackrel \pi \arrow S$
its twistor space, and $I, L\in S$ 
complex structures which satisfy 
$V^{0,1}_L\cap V^{0,1}_I=0$. Consider the category
$\c C_I$ constructed above, and let
$\c C_I\stackrel\Phi\arrow \Coh(T,L)$ be the 
restriction map. Then $\Phi$ is equivalence 
of categories.

\hfill

{\bf Proof:} We construct the inverse functor 
$\Upsilon:\;\c C_I\arrow  \Coh(T,L)$
as follows. Take $F_L\in \Coh(T,L)$.
Let $B_L$ be its reflexization, which
is smooth by \ref{_gene_refle_Proposition_},
and $B\in \c C_I$ the corresponding
bundle over the twistor space,
which is defined in a canonical 
way by \ref{_restri_from_Fl_I_equi_Theorem_}.

Let $x_0\in T$ be a singular point of $F_L$.
Consider the infinitesimal neighbourhood
of 
\[ \{x_0\}\times(S\backslash I)\subset \Tw(T) \backslash (T,I)
\] 
in $\Tw(T) \backslash (T,I)$, denoted, as in Subsection
\ref{_local_trivi_Subsection_}, by $\Tw(T)_{x_0, I}$.

Then $B$ is trivialized over $\Tw(T)_{x_0, I}$,
and we may write $B\restrict{\Tw(T)_{x_0, I}} = \xi^*(B_{L, x_0})$,
where $B_{L, x_0}$ is the germ of $E_L$ in $x_0$.
Let $F\restrict{\Tw(T)_{x_0, I}}:=\xi^*(F_{L, x_0})$ be 
the sheaf corresponding to $F_L$. 
Since $F_L$ has isolated singularities, the sheaves
$F\restrict{\Tw(T)_{x_0, I}}$ and
$B\restrict{\Tw(T)_{x_0, I}}$ are 
canonically isomorphic outside of
$\{x_0\}\times(S\backslash I)$.
Gluing together $F\restrict{\Tw(T)_{x_0, I}}$
and $B$, we obtain a sheaf $F_0$ on $\Tw(T) \backslash (T,I)$.
Let $A$ be the singular set of $F_L$. Then
$B$ is equal to $F_0$ outside of 
$\{A\} \times S\subset \Tw(T)$,
hence $F_0$ can be extended smoothly
{}from $\Tw(T) \backslash (T,I)$
to  a sheaf $F_1$ on $\Tw(T) \backslash Z$, where
$Z= \{A\} \times \{I\}$ is the finite set
defined in \ref{_C(I)_Definition_}.
Consider the sheaf $F:=j_*(F_1)$, where 
$j:\; \Tw(T) \backslash Z \hookrightarrow \Tw(T)$
is the natural embedding. The correspondence
$F_L \arrow F$ is clearly functorial.

The sheaf $F$ by construction 
belongs to $\c C_I$, and satisfies
$\Phi(F)=F_L$. Therefore, $F_L\arrow F$
gives an inverse functor of 
$\c C_I\stackrel\Phi\arrow \Coh(T,L)$.
We proved \ref{_C(I)_equiv_Theorem_}.
\endproof

\hfill

The following corollary immediately
follows from \ref{_C(I)_equiv_Theorem_}.

\hfill

\corollary
Let $T$ be a generic compact torus, equipped with
flat Hermitian metric $g$, and $S$ the set
of complex structures compatible with $g$.
Then, for any generic complex structures
$I, J\in S$, the categories of coherent sheaves
$\Coh(T,I)$ and $\Coh(T,J)$ are isomorphic.

\endproof


\section{Moduli of compact tori}
\label{_moduli_conne_Section_}

Let $T= \C^n/\Z^{2n}$ be a compact torus
with a flat Riemannian metric $g$. We have shown that
for all generic complex structures $I, J$
compatible with $g$, the category
$\Coh(T,I)$ is equivalent to $\Coh(T,J)$.
To prove that $\Coh(T,I)$ is equivalent to $\Coh(T,J)$
for arbitrary generic complex structures, 
it suffices to prove the following proposition.

\hfill

\proposition\label{_conne_moduli_6_hops__Proposition_}
Let $T$ be a compact torus, $T=V/\Z^{2n}$, $V\cong \C^n$,
and $R$ its marked moduli space (Section \ref{_gene_tori_Section_}).
Given $I\in R$, pick a flat Riemannian metric compatible
with $I$, and let $R_g\subset R$ be the set of
all complex structures $J\in R$ compatible
with $g$. Take a generic complex structure
$I_1 \in R_g$, and repeat the same procedure,
obtaining $I_2$. Then, after at most 6 iterations,
we can obtain any generic complex structure
on $T$ from $I$.

\hfill

{\bf Proof:} We use the following Claim.
Please note that the flat metrics on $T$ are
in one-to-one correspondence with the
Euclidean metrics on $V$.

\hfill

\claim\label{_goth_G_Claim_}
In assumptions of 
\ref{_conne_moduli_6_hops__Proposition_},
let $\c G$ be the set of all Euclidean metrics $g_1$ on $V$
such that $g_1$ is compatible with some complex structure
$I\in R_g$. Then $g_1\in \c G$ if and only if all
eigenvalues of the symmetric matrix $g_1 g^{-1}$
occur in pairs. In other words, $g_1 \in \c G$ if
and only if there is an orthonormal basis in $(V,g)$
such that $g_1$ is written in this basis as
\begin{equation}\label{_eigenva_in_pairs_Equation_}
\left(
\begin{matrix}
\alpha_1 &&&&&&\\
&\alpha_1&&&&&\\
&&\alpha_2&&&&\\
&&&\alpha_2&&&\\
&&&&\ddots &&\\
&&&&&\alpha_n&\\
&&&&&&\alpha_n
\end{matrix}
\right)
\end{equation}
{\bf Proof:}
``If'' part. Let $I_1\in R_g$ be a complex structure
such that $g_1$ is compatible with $I_1$. Then 
$g, g_1$ are Hermitian forms on the complex
vector space $(V,I_1)$. A standard linear-algebraic
argument implies
that the corresponding Hermitian metrics have a common 
orthogonal basis $z_1, z_2, ..., z_n\in (V,I_1)$.
After rescaling, we may assume that $z_1, z_2, ..., z_n$
is orthonormal with respect to $g$. Consider the
corresponding basis $z_1, I_1(z_1), z_1, I_1(z_1), ...$
in $V$ over $\R$. In this
basis, $g$ is written as identity matrix,
and $g_1$ as \eqref{_eigenva_in_pairs_Equation_}.
This proves the ``if'' part of \ref{_goth_G_Claim_}.

Conversely, assume that $g^{-1}g_1$ has eigenvalues occuring
in pairs. Finding a common orthogonal basis, which is
orthonormal in $g$, we find that $g_1$ is written in this basis
as \eqref{_eigenva_in_pairs_Equation_}. Let $I_1$ act in this
basis as
\begin{equation}
\left(
\begin{matrix}
 0&1&&&&&\\
-1&0&&&&&\\
   && 0&1&&&\\
   &&-1&0&&&\\
&&&&\ddots &&\\
&&&&&0&1\\
&&&&&-1&0
\end{matrix}
\right)
\end{equation}
Then $I_1$ is a complex structure compatible with $g$ and $g_1$.
Therefore, $g_1\in \c G$. This proves \ref{_goth_G_Claim_}.
\endproof

\hfill

The following elementary claim is proven in the same
way as one proves \ref{_non-generic_Proposition_}.

\hfill

\claim\label{_in_R_g_all_gene_Claim_}
In assumptions of \ref{_goth_G_Claim_}, let 
$R^\circ_g$ be the set of all generic complex structures
$I_1\in R_g$. Then $R_g \backslash R^\circ_g$
has measure zero.

\endproof

\hfill

Using \ref{_in_R_g_all_gene_Claim_}, 
\ref{_goth_G_Claim_} and Fubini theorem, we obtain

\hfill

\corollary\label{_c_G^circ_almost_all_Corollary_}
In the above assumptions, let $\c G^\circ$ 
be the set of all Euclidean metrics $g_1\in \c G$
such that $g_1$ is compatible with $I_1\in R^\circ_g$.
Then $\c G\backslash \c G^\circ$ has measure
zero. 

\endproof

Now we can prove the following lemma.

\hfill

\lemma\label{_after_3_ite_almost_all_Lemma_}
In the assumprions of \ref{_conne_moduli_6_hops__Proposition_},
let $R_3(I)$ be the set of all $I_3\in R$ 
which can be obtained after three iterations
of the procedure defined in
\ref{_conne_moduli_6_hops__Proposition_}.
Then $R \backslash R_3(I)$ 
has measure zero.

\hfill

{\bf Proof:} Let $\c F\in \End V$
be the set of all symmetric operators on $(V, g)$
which can be written as $U A U^{-1}$,
where $U$ is orthogonal, and $A$ takes form 
\eqref{_eigenva_in_pairs_Equation_}. From
\ref{_goth_G_Claim_} it is clear that $\c G = \c F\cdot g$.
The second iteration of this procedure gives the following
trivial claim.

\hfill

\claim
In the above notation, let $\c G_1$ be the set of all
Euclidean metrics on $V$ which are compatible
with some $I_2\in \R$, with $I_2$ being
compatible with some $g_1\in \c G$. Then
$\c G_1 = \c F \cdot \c F\cdot g$.

\endproof

\hfill

Clearly, every Euclidean metric can be obtained this way; even if
$A_1$ and $A_2$ commute, $U_1 A_1 A_2 U_1^{-1}$
runs through the set of all symmetric matrices.  
Indeed, any diagonal matrix can be obtained as a product
of two diagonal martices with eigenvalues occuring in pairs.

Let $\c G_1^\circ$ be the set of all metrics
compatible with a generic $I_2\in R$, after
a second iteration of the procedure defined
in \ref{_conne_moduli_6_hops__Proposition_}.
Then $\c G_1 \backslash G_1^\circ$ has measure 0,
as one can easily surmise from
\ref{_c_G^circ_almost_all_Corollary_}.
Then, $\c G_1$ is the set of all 
Euclidean metrics except, possibly, 
 a measure zero set, and $\c G_1$
is compatible with all $I\in R$ except,
possibly, a measure zero set.
This proves \ref{_after_3_ite_almost_all_Lemma_}.
\endproof

\hfill

Return to the proof of \ref{_conne_moduli_6_hops__Proposition_}.
Let $I, J\in R$ be generic complex structures,
and $R_3(I)$, $R_3(J)\subset R$ the sets associated
with $I$, $J$ as in \ref{_after_3_ite_almost_all_Lemma_}.
The sets $R\backslash R_3(I)$, $R\backslash R_3(J)$
have measure zeto, as \ref{_after_3_ite_almost_all_Lemma_}
implies. Therefore, $R_3(I)$ and $R_3(J)$ have
non-empty intersection. Take a chain
\[ I \to I_1 \to I_2 \to I_3 = J_3 \to J_2 \to J_1\to J
\]
of generic complex structures with each successive
pair having a Euclidean metric compatible with both.
We obtain that $J$ can be obtained from $I$ after
6 iterations of the procedure defined in 
\ref{_conne_moduli_6_hops__Proposition_}.
\ref{_conne_moduli_6_hops__Proposition_} is proven.
\endproof

\hfill

\hfill

{\bf Acknowledgements:} I am grateful to 
F. Bogomolov for valuable advice
and D. Kaledin for an important insight in the 
geometry of $SO(2n)/U(n)$.

\hfill

{\small

\noindent {\sc Misha Verbitsky\\
University of Glasgow, Department of Mathematics, 15
  University Gardens, Glasgow, Scotland.}\\
\tt verbit@maths.gla.ac.uk, \ \  verbit@mccme.ru 

}

\end{document}